\documentclass{amsart}
\usepackage{amsthm}
\usepackage{graphicx,textcomp}
\usepackage{amsmath,bm,amsfonts,mathrsfs,amssymb}
\usepackage{hyperref}

\usepackage{appendix}
\usepackage[normalem]{ulem}
\usepackage{graphicx,color}
\usepackage{enumitem}   
\usepackage{xcolor}

\begin{document}
\title{Variational formulas for determinant of Laplacian on higher genus polyhedral surface}
\author{Dmitrii Korikov}
\email{dmitrii.v.korikov@gmail.com}
\author{Alexey Kokotov}
\email{alexey.kokotov@concordia.ca}
	

\maketitle
\begin{abstract}
Let $X$ be a Riemann surface of genus $g\ge 1$ endowed with a flat conical metric $m$ and let ${\rm det}\,\Delta$ be the $\zeta$-regularized determinant  of the Friedrichs Laplacian on $(X,m)$. We derive variational formulas for ${\rm det}\,\Delta$ with respect to conical points and conical angles within a given conformal class. Integration of them leads to an explicit expression for ${\rm det}\,\Delta$ up to moduli dependent factor. The latter, in principle, can be calculated via comparison of the above result with the well-known formulas for the case of flat conical metrics with trivial holonomy.
\end{abstract}

Let $X$ be a Riemann surface of genus $g\ge 1$ and $m$ be a conformal flat conical metric on $X$. Let $\Delta$ be a Friedrichs Laplacian on $X$ corresponding to the metric $m$. We study the dependence of the zeta-regularized determinant ${\rm det}\Delta$ of $\Delta$ on the positions of the conical singularities and conical angles.

For higher genus case, the explicit formulas for ${\rm det}\,\Delta$ were obtained in \cite{KoKorotTau} for metrics $m$ with trivial holonomy. The latter metrics are of the form $m=|\omega|^2$, where $\omega$ is an Abelian differential on $X$; in particular, all the conical angles are integer multiples of $2\pi$. For general polyhedral metrics, the Polyakov type comparison formula 
\begin{equation}
\label{KokCHS fla}
\frac{{\rm det}\Delta^{(m)}/{\rm Area}(X,m)}{{\rm det}\Delta^{(\tilde{m})}/{\rm Area}(X,\tilde{m})}=\frac{\prod_{k=1}^M C(\beta_k)\prod_{j=1}^{\tilde{M}}f_j}{\prod_{j=1}^{\tilde{M}} C(\tilde{\beta}_j)\prod_{j=1}^{\tilde{M}}\tilde{f}_k}
\end{equation}
for the determinants of the Laplacians $\Delta^{(m)},\Delta^{(\tilde{m})}$ corresponding to two conformally equivalent flat conical metrics $m,\tilde{m}$, respectively, was obtained in \cite{KokPAMS}. Here $\beta_k$, $\tilde{\beta}_j$ are the angles at the conical singularities $P_k$,$\tilde{P}_k$ of the metrics $m$, $\tilde{m}$, respectively, while the coefficients $f_j$ (resp., $\tilde{f}_j$) are calculated by representation of the one metric $m$ (resp., $\tilde{m}$) near $\tilde{P}_j$ (resp., $P_k$) in the distinguished coordinates corresponding to the second metric $\tilde{m}$ (resp., $m$). The $C(\beta_k),C(\tilde{\beta}_k)$ admit explicit expressions in terms of the determinants of the Dirichlet Laplacians on the finite cones with conical angles $\beta_k,\tilde{\beta}_k$; the explicit formulas for the latter are obtained by Spreafico \cite{Spreafico} via involved calculations with special functions.

The derivation of (\ref{KokCHS fla}) is based on the application of the classical Alvarez comparison formula and the BFK gluing formula for determinants \cite{BFK}. Formula (\ref{KokCHS fla}) can be considered as an analogue of the global Polyakov formula for the case of flat conical metrics. The the classical global Polyakov formula is obtained via the integration of the infinitesimal one that describes the variations of the determinant under the infinitesimal variation of (smooth) conformal factor. 

The goal of this paper is to derive an analogue of the infinitesimal
Polyakov formula for the case of flat conical metrics and to derive expressions for the determinant avoiding the use of the BFK formula and complicated calculations with special functions. In the genus 0 case, this was done in \cite{KoKoPolyPoly}. 

In the present paper, we generalize the argument of \cite{KoKoPolyPoly} for the surfaces of arbitrary genera $g>0$. First, in Section \ref{sec eig}, we study the dependence of individual eigenvalues of $\Delta$ on positions of the singularities and the conical angles; to this end, the classical perturbation theory for elliptic equation on manifolds with singularities (provided by \cite{MNP}) is applied. As a result, we prove variational formulas for the zeta function $\zeta_{\Delta-\lambda}(2)$ of the operator $\Delta-\lambda$ ($\lambda\in\mathbb{C}$). The asymptotics of the latter as $\Re\lambda\to-\infty$ is derived (in Section \ref{heat kernel sec}) by careful analysis of the asymptotics of the heat kernel near conical singularities. Then one can apply the inverse Mellin transform to obtain variational formula (\ref{zeta calc}) for ${\rm det}\Delta$ (see Section \ref{sec det}). The coefficients in (\ref{zeta calc}) are calculated in Section \ref{sec coeff} for the case $g>1$ and in Section \ref{sec tori} for the case $g=1$ by means of explicit expressions (\ref{explicit metric}) for the polyhedral metric found in \cite{KoKo-poly} (see also \cite{KokLN,KokPAMS}). Finally, the integration of (\ref{zeta calc}) leads to explicit expression (\ref{Main fomula 2}) for the determinant ${\rm det}\Delta$ with the ingredients provided by (\ref{Main formula 1}), (\ref{scaling formula}), and (\ref{beta part}), (\ref{beta part torus}). For the higher genus case, the only undetermined constant ${\bf c}_0$ in (\ref{Main fomula 2}) depends on the moduli of $X$ and could be, in principle, calculated by substitution the metric with trivial holonomy into (\ref{Main fomula 2}) and then comparing the result with formulas obtained in \cite{KoKorotTau}; this last step is omitted here. In the case $g=1$, explicit value (\ref{torical const}) of ${\bf c}_0$ is easily calculated by comparison (\ref{Main fomula 2}) with the Ray-Singer formula \cite{RS}.

We hope that our approach will work also for Laplacians acting in hermitian linear bundles (say, powers of $K$).

\section{Family of polyhedral metrics.} Let $X$ be a genus $g$ Riemann surface and let $\mathcal{D}:=\sum_{k=1}^M b_k P_k$ be the generalized (i.e., $b_k$ are not necessarily integer) divisor on $X$ obeying $\sum_{j=1}^M b_j=2g-2$. Due to the Troyanov theorem \cite{Troy}, there is the unique (up to homothety) conformal metric $m$ on $X$ obeying the following conditions
\begin{enumerate}
\item $m$ is flat on $\dot{X}:=X\backslash\{P_1,\dots,P_M\}$. Representing the metric in holomorphic local coordinates,
\begin{equation}
\label{metric}
m=e^{-\phi(z)}|dz|^2,
\end{equation}
one obtains that $\phi$ is harmonic, $\partial_z\partial_{\overline{z}}\phi=e^{\phi(z)}K(z)/2=0$ (where $K$ is the Gaussian curvature of $m$).

\item In the neighborhood of each $P_j$, the metric admits the representation 
\begin{equation}
\label{metric near vertices}
m=|e^{-u_j(z)}|\cdot|(z-z(P_j))^{b_j}dz|^2
\end{equation}
(in any holomorphic local coordinates), where $u_j$ is smooth. Due to condition 1., $\Re u_j(z)=\phi(z)+2b_j{\rm log}|z-z(P_j)|$ is harmonic and, therefore, one can chose $u_j$ to be holomorphic. Then a neighborhood of $P_j$ is isometric to a neighborhood of the vertex of the cone $\mathbb{K}_\beta$ of opening $\beta=\beta_j:=2\pi(1+b_j)$ by means of the map
\begin{equation}
\label{isometry to cone}
P\mapsto\xi_j(P):=\int_{z(P_j)}^{z(P)}(z-z(P_j))^{b_j(t)}e^{-u_j(z)/2}dz.
\end{equation}
\end{enumerate}
The metric $m$ obeying the above conditions for some $\mathbb{D}$ is called {\it polyhedral}. 

From now on we allow $m$ to depend on a (real or complex) parameter $\tau$. We assume that the coefficients $b_k$, $z(P_k)$ of the divisor $\mathbb{D}$ depend smoothly on $\tau$. Due to explicit formula (\ref{explicit metric}) for polyhedral metrics below, this implies that the functions $\phi|_{\dot{X}}$ and $u_j$ in the above conditions are smooth in $z,\tau$. In the subsequent, as a rule, we omit the dependence on $\tau$; the derivative $\partial_\tau$ will be denoted by a dot. 

Note that $\dot{\phi}(z(P))$ is independent of the choice of the local coordinates $z$ of $P$ and is harmonic on $\dot{X}$ due to condition 1 while (\ref{isometry to cone}) implies
\begin{align*}
\frac{1}{z-z(P_j)}=e^{-u_{j}(z(P_j))\pi /\beta_j}\Big(\frac{\beta_j\xi_j}{2\pi}\Big)^{-\frac{2\pi}{\beta_j}}\,-\frac{u'_{j}(z(P_j))}{2(b_j+2)}+o(1).
\end{align*}
Then $\dot{\phi}(z)=\partial_\tau[\Re u_j(z)-2b_j{\rm log}|z-z(P_j)|]$ admits the expansion
\begin{equation}
\label{Fourier}
\dot{\phi}=A_j(\tau)+\tilde{A}_j(\tau){\rm log}|\xi_j|+\sum_{k=-1,1,2,\dots}\big[c_{k,j}(\tau)\xi_j^{\frac{2\pi k}{\beta_j}}+c^\dag_{k,j}(\tau)\overline{\xi}_j^{\frac{2\pi k}{\beta_j}}\big],
\end{equation}
near $P_j$, where
\begin{align}
\label{Fourier coef}
\begin{split}
\tilde{A}_j=-\frac{2\dot{\beta}_j}{\beta_j}, \quad A_j&=-\frac{b_j}{2(b_j+2)}\Big[u'_j(z(P_j))\frac{\partial z(P_j)}{\partial \tau}+\overline{u'_j(z(P_j))}\frac{\partial \overline{z(P_j)}}{\partial \tau}\Big]-\\
&-\frac{\dot{\beta}_j}{\beta_j}\Big[2{\rm log}\Big(\frac{\beta_j}{2\pi}\Big)+\Re u_j(z(P_j))\Big]+[\partial_\tau\Re u_j](z(P_j)).
\end{split}
\end{align}

\subsubsection{Explicit formulas for polyhedral metrics}
The explicit expressions for the polyhedral metrics on surfaces of genera $g\ge 2$ are proposed in \cite{KokLN,KokPAMS} (under the assumption that none of $P_j$ is a Weierstrass point) and in Section 2, \cite{KoKo-poly} (for the general case). According to the latter result, any polyhedral metric (\ref{metric}) on $X$ can be written as
\begin{align}
\label{explicit metric}
\begin{split}
\phi(z)=-\sum_{k=1}^M b_k\Big[{\rm log}\Big|\frac{E(z,z'_k)}{E(p_1,z'_k)}\Big|^2+Q(z,z'_k)\Big]-{\rm log}C-4{\rm log}|\sigma(z,p_0)|.
\end{split}
\end{align}
Here $\mathcal{D}=\sum_{k=1}^M b_k P_k$ is the generalized divisor of the metric,
\begin{align}
\label{Q form}
Q(z,z')=\frac{4\pi}{(g-1)^2}(\Im K^{z})^T(\Im\mathbb{B})^{-1}\Im K^{z'};
\end{align}
each $z_k$ is arbitrary holomorphic coordinate near the vertex $P_k$ while $z'_k=z_k(P_k)$, $p_0,p_1$ are arbitrarily fixed points, $p_1$ does not coincide with vertices, $C>0$ is the scaling factor, $E$ is the prime form, $\mathscr{A}$ and $K^z$ is the Abel transform and the vector of Riemann constants, respectively, that correspond to the basis $\vec{v}=(v_1,\dots,v_g)^T$ of Abelian differentials and the basepoint $z$, and the function $\sigma$ is well-defined on $X$ by the rule
\begin{equation}
\label{Troy 2}
\sigma(z,p_0)=\frac{\theta\Big(\mathscr{A}\big(\sum_{j=1}^g x_j-gz\big)+K^z\Big)}{\theta\Big(\mathscr{A}\big(\sum_{j=1}^g x_j-p_0-(g-1)z\big)+K^z\Big)}\prod_{j=1}^g\frac{E(x_j,p_0)}{E(x_j,z)}
\end{equation}
with arbitrary $x_j$ (see (1.13), \cite{Fay}). Formula (\ref{explicit metric}) is obtained by substitution of formulas (18) and (22), \cite{KoKo-poly} into (\ref{metric}), taking into account the identity 
$$K^z=(g-1)\mathscr{A}(z-z')+K^{z'}$$
from Proposition 1.1, \cite{Fay} and adding manually some additional (independent of $z$) terms to make the metric independent of the choice of the coordinates $z'_k$ of $P_k$.

Comparing(\ref{explicit metric}) with (\ref{metric near vertices}) yields
\begin{align}
\label{u j explict}
\begin{split}
u_j(z_j)=&-\sum_{k=1}^M b_k\Bigg[Q(z,z'_k)+{\rm log}\Big[\frac{E(z_j,z'_k)}{E(p_1,z'_k)}\Big]^2\Bigg]+\\
&+2b_j{\rm log}(z_j-z'_j)-{\rm log}C-4{\rm log}\sigma(z_j,p_0)
\end{split}
\end{align}

\section{Variation of individual eigenvalues and $\zeta_{\Delta-\mu}(2)$}
\label{sec eig}
Family of metrics (\ref{metric}) defines the family of the Friedrichs Laplacians on $X$
$$\tau\mapsto\Delta_\tau=e^{\phi_\tau(z)}\partial_z\partial_{\overline{z}}.$$
Suppose that $\lambda_{\tau_0}$ is an eigenvalue of $\Delta_{\tau_0}$ of multiplicity $d$. Then it can be embedded into the eigenpair families $\tau\mapsto (\lambda^{(k)}_{\tau},u^{(k)}_{\tau})$, $k=1,\dots,d$ obeying
\begin{enumerate}
\item $\tau\mapsto\lambda^{(k)}_{\tau}$ is smooth and $(z,\tau)\mapsto u^{(k)}_{\tau}(z)$ is smooth outside vertices $\tau\mapsto P_j(\tau)$ (here $z$ is an arbitrary holomorphic coordinate on $SX$);
\item the (convergent) asymptotics series
$$u^{(k)}_{\tau}(\xi_j)=\sum_{s,m=0}^{\infty}[c_{j,s,l}(\tau)\xi_j^{\frac{2\pi s}{\beta}+l}\overline{\xi}_j^l+c^\dag_{j,s,m}(\tau)\xi_j^l\overline{z}^{\frac{2\pi s}{\beta}+l}]$$
for $u^{(k)}_{\tau}$ near $P_j$ admit the differentiation with respect to $t$;
\item $u^{(k)}_{\tau_0}$, $k=1,\dots,d$ constitute the orthonormal basis in ${\rm Ker}(\Delta_{\tau_0}-\lambda_{\tau_0})$.
\end{enumerate}
The above property can by proven in the following way (see Lemma 4.1 and Appendix B, \cite{KoKoPolyPoly}). By cutting off the metric $\epsilon$-neighborhoods $\mathbb{K}_j(\epsilon)$ of conical points $P_j$ and applying the explicit formulas for the DN maps of $\mathbb{K}_j(\epsilon)$, one reduces the original eigenproblem to the problem in the domain $X^\circ_\tau(\epsilon)$; the latter can be considered as a fixed domain endowed with the smooth families of metric and (non-local) boundary operators. For the latter problem, one applies the standard perturbation theory to obtain 1. Hence, by the use of the explicit formulas for the resolvent kernels in $\mathbb{K}_j(\epsilon)$, one obtains 2. It worth noting that the same statement can be obtained from the general theory of elliptic problems in singularly perturbed domains (see Chapters 4 and 6, \cite{MNP}). 

Let $\tau\mapsto(\lambda_\tau,u_\tau)$ be one of the above eigenpair families. Due to property 1., one can differentiate the equation $(\Delta-\lambda)u=0$ with respect to $\tau$ to obtain
\begin{equation}
\label{variation of eigenfucntions}
(\Delta-\lambda)\dot{u}=(\dot{\lambda}k-\dot{\Delta})u=(\dot{\lambda}-\dot{\phi}\Delta)u=(\dot{\lambda}-\lambda\dot{\phi})u.
\end{equation}
Differentiating the asymptotics of $u_\tau$ near vertices with respect to $\tau$, one conclude that right-hand side of the Green formula
\begin{align*}
((\Delta-\lambda)\dot{u},u)_{L_2(X^\circ_\tau(\epsilon);m_\tau)}-(\dot{u},(\Delta-\lambda)u)_{L_2(X^\circ_\tau(\epsilon);m_\tau)}=\\
=(\dot{u},\partial_{\nu_\tau}u)_{L_2(\partial X^\circ_\tau(\epsilon);m_\tau)}-(\partial_{\nu_\tau}\dot{u},u)_{L_2(\partial X^\circ_\tau(\epsilon);m_\tau)}
\end{align*} 
tends to zero as $\epsilon\to 0$. Now the substitution of (\ref{variation of eigenfucntions}) into the Green formula, passing to the limit as $\epsilon\to 0$ and applying property 3. yield
\begin{equation}
\label{variation of eigenvalues}
\frac{\dot{\lambda}}{\lambda}={\rm p.v.}\int_X\dot{\phi}u^2dS,
\end{equation}
where $dS=e^{-\phi}d\overline{z}\wedge dz/2i$ is the area element. Here the Cauchy principal value should be taken since the singularities $c_{-1,j}\xi^{-\frac{2\pi}{\beta_j}}+c^\dag_{-1,j}\overline{\xi}^{-\frac{2\pi}{\beta_j}}$ of $\dot{\phi}$ near vertices $P_j$ (see (\ref{Fourier})) are not absolutely integrable. 

As a corollary of (\ref{variation of eigenvalues}), we have
$$\partial_\tau((\lambda-\mu)^{-2})={\rm p.v.}\int\limits_X \partial_\mu^2\Big(\mu \frac{u^2}{\mu-\lambda}\Big)\dot{\phi}dS$$
for any $\mu\in\mathbb{C}$, $\mu\ne\lambda_\tau$. Making the summation over all eigenvalues of $\Delta=\Delta_\tau$ in the last formula, one arrives at
$$\dot{\zeta}_{\Delta-\mu}(2)={\rm p.v.}\int\limits_X\partial_\mu^2[-\mu R_\mu(x,y)]\Big|_{y=x}\dot{\phi}(x)dS(x),$$
where $R_\mu(x,y)=R_{\mu,\tau}(x,y)$ is the resolvent kernel of $\Delta=\Delta_\tau$. Here one can add any linear function of $\mu$ to the expression in the square brackets; in particular, it can be replaced by
\begin{align*}
\mu \Big[\frac{1}{2\pi}K_0(d\sqrt{-\mu})-R(x,y)-\frac{K_0(d\sqrt{-\mu})-{\rm log}d}{2\pi}\Big]-\frac{1}{A},
\end{align*}
where $A$ is the area of $(X_t,m_t)$, $d=d(x,y)$ is the distance between $x,y\in X$ in the metric $m_t$, and $K_0$ is the Macdonald function. The diagonal value of the last expression is equal to $\psi_\mu(x)+\tilde{\psi}_\mu$, where
\begin{align}
\label{psi mu}
\psi_\mu(x)=\mu \Big[\frac{K_0(d\sqrt{-\mu})}{2\pi}-R_\mu(x,y)\Big]_{y=x}, \ \tilde{\psi}_\mu=\frac{\mu\,({\rm log}(4|\mu|)+2\gamma)}{4\pi}-\frac{1}{A}.
\end{align}
Thus, one arrives at
\begin{equation}
\label{zeta of 2}
\dot{\zeta}_{\Delta-\mu}(2)=\partial_\mu^2\Big[{\rm p.v.}\int\limits_X[\psi_\mu+\tilde{\psi}_\mu]\dot{\phi}dS\Big].
\end{equation}
Note that $\psi_\mu(x)+\tilde{\psi}_\mu\to 0$ as $\mu\to 0$ uniformly in $x$ due to the standard expansion $R_\mu=-(1/\mu A)+G+O(\mu)$ of the resolvent kernel (where $G$ is the Green function of $\Delta$).

\section{Variation of ${\rm det}\Delta$}
\label{sec det}
The zeta functions $\zeta_{\Delta-\mu}(2)$ and $\zeta_{\Delta}(s)$ are connected via
\begin{equation}
\label{zeta of s via zeta of 2}
\zeta_{\Delta}(s)=\frac{\mathcal{M}^{-1}\{\zeta_{\Delta-\mu}(2)\}(s-1)}{s-1},
\end{equation}
where $\mathcal{M}^{-1}$ is the inverse Mellin transform 
$$\mathcal{M}^{-1}\{f(\mu)\}(s):=\frac{1}{2\pi i}\int_\Gamma \mu^{-s}f(\mu)d\mu$$ 
(here $\Gamma $ is the contour enclosing the cut $(-\infty,0)$).

Now we make use the following elementary property of the Mellin transform (see Lemma 5.1, \cite{KokKorWZ} and Lemma 3.1 and Appendix A, \cite{KoKoPolyPoly}). Suppose that $f$ is holomorphic near $(-\infty,0)$ and obeys the (admitting differentiation) asymptotics
\begin{equation}
\label{F asymp}
f(\mu)=\sum_{k=1}^K [f_k+\tilde{f}_k\mu{\rm log}(-\mu)]\mu^{r_k}+O(\mu^{\kappa}), \quad \Re\mu\to -\infty
\end{equation}
with $\kappa<0$, $r_k\in\mathbb{R}$. Then the function 
$$\hat{f}(s):=\frac{\mathcal{M}^{-1}\{f(\mu)\}(s-1)}{s-1}$$
obeys
$$\hat{f}(0)=-[f]^{{\rm log}}, \quad -\partial_s\hat{f}(0)=[f]^{\infty}-f(0),$$
where $[f]^{\infty}$ and $[f]^{{\rm log}}$ denote the constant and logarithmic term in (\ref{F asymp}), respectively.

Now we assume that $f(\mu)=\dot{\zeta}_{\Delta-\mu}(2)$ admits asymptotics (\ref{F asymp}); this assumption will be justified by formulas (\ref{integrated asymptotic const term}) and (\ref{log integral at infty}) below. Differentiating (\ref{zeta of s via zeta of 2}) with respect to $\tau$, taking into account (\ref{zeta of 2}), and applying the above property, one arrives at
\begin{equation}
\label{zeta calc 1}
-\partial_s\dot{\zeta}_{\Delta}(0)=\Big[{\rm p.v.}\int\limits_X[\psi_\mu+\tilde{\psi}_\mu]\dot{\phi}dS\Big]^{\infty}=\Big[{\rm p.v.}\int\limits_X\psi_\mu\dot{\phi}dS\Big]^{\infty}+\frac{\dot{A}}{A}\,.
\end{equation}
Here the last equality follows from (\ref{psi mu}) and the equality $\dot{A}=-\int_X\dot{\phi}dS$. 

Far from vertices, the parametrix for the resolvent kernel as $\Re\mu\to-\infty$ is provided by the Macdonald function, i.e.,
$$\frac{K_0(d(x,y)\sqrt{-\mu}))}{2\pi}-R_\mu(x,y)=O(\mu^{-\infty}), \qquad \Re\mu\to-\infty$$
uniformly in $x,y$ separated from vertices. As a corollary, $\psi_\mu(x)$ exponentially decays as $\Re\mu\to-\infty$ uniformly in $x$ separated from  vertices. Thus, one can replace $X$ in the last integral in (\ref{zeta calc 1}) with the union of small $\epsilon$-neighborhoods $\mathbb{K}_j(\epsilon)$ of conical points $P_j$. As a result, one obtains
\begin{equation}
\label{zeta calc 2}
\partial_\tau[{\rm log}({\rm det}\Delta/A)]=\sum_{j=1}^M \Big[{\rm p.v.}\int_{\mathbb{K}_j(\epsilon)}\psi_\mu\dot{\phi}dS\Big]^{\infty}.
\end{equation}

To derive the asymptotics the integrals in the right-hand side of (\ref{zeta calc 2}), one needs the asymptotics of the resolvent kernel $R_\mu(x,y)$ as $\Re\mu\to-\infty$ which is uniform in $x,y$ arbitrarily close to the vertex $P_j$. To this end, recall that the heat kernel for the (Friedrichs) Laplacian in the infinite cone $\mathbb{K}_\beta$ with opening angle $\beta$ is given by
\begin{equation}
\label{HK Cone}
{\rm H}_{t,\beta}(x,x')=\frac{1}{8\pi i\beta t}\int\limits_{\mathscr{C}}{\rm exp}\Big(-\frac{\mathfrak{r}^2(r,r',\vartheta)}{4t}\Big)\,\Xi_\beta(\vartheta,\varphi,\varphi') d\vartheta
\end{equation}
(see \cite{Carslaw,Dowker}). Here $(r,\varphi)$ and $(r',\varphi')$ are the polar coordinates of the points $x$ and $x'$, respectively, and
\begin{equation*}
\mathfrak{r}^2:=r^2-2rr'{\rm cos}\vartheta+r'^{2}, \qquad \Xi:={\rm cot}\Big(\pi\beta^{-1}(\vartheta+\varphi-\varphi')\Big).
\end{equation*}
The contour $\mathscr{C}$ in (\ref{HK Cone}) belongs to the strip $\Re\theta\in(-\pi,\pi)$ and consists of two lines $\pm l:=[\pm\pi\mp i\infty,\pm\pi\pm i\infty]$ and arbitrarily small anti-clockwise circles centered at the roots of $\Xi$. The resolvent kernel ${\rm R}_{\mu,\beta}(x,x')$ for the Laplacian in $\mathbb{K}_\beta$ is obtained by the Laplace transform in (\ref{HK Cone}), ${\rm R}_{\mu,\beta}=\int_{0}^{+\infty}e^{\mu t}{\rm H}_{t,\beta}dt$.

Let us represent the resolvent kernel $R_\mu(x,y)$ as 
\begin{equation}
\label{parametrix resolvent}
R_\mu(x,x')=\chi_j(x){\rm R}_{\mu,\beta_j}(\xi_j(x),\xi_j(x'))+\tilde{R}_\mu(x,x'),
\end{equation}
where $\chi_j$ is the cut-off function equal to one in $\mathbb{K}_j(2\epsilon)$ and to zero outside $\mathbb{K}_j(3\epsilon)$. It easily follows from (\ref{HK Cone}) that $\tilde{R}_\mu(\cdot,x')\in {\rm Dom}\Delta$. In view of (\ref{HK Cone}), ${\rm H}_{t,\beta}(x,x')$ decays exponentially as $t\to +0$ uniformly in $x$ separated from $x'$. Hence, due to the properties of the Laplace transform, ${\rm R}_{\mu,\beta}(x,x')=O(|\mu|^{-\infty})$ as $\Re\mu\to-\infty$ uniformly in $x$ separated from $x'$. If $x'\in \mathbb{K}_j(\epsilon)$, then the support of the function $(\Delta-\mu)\tilde{R}_\mu(\cdot,x')=[\chi_j,\Delta]R^{(0)}_\mu(\cdot,x')$ is separated from $x'$ and, thus, 
\begin{equation}
\label{parametrix rhs estimate}
(\Delta-\mu)\tilde{R}_\mu(\cdot,x')=O(|\mu|^{-\infty}), \qquad \Re\mu\to-\infty
\end{equation}
uniformly in $x'\in\mathbb{K}_j(\epsilon)$ and $x\in X$; the same estimate is valid for all derivatives of $(\Delta-\mu)\tilde{R}(\cdot,x')$. In view of the standard operator estimate $(\Delta-\mu)^{-1}=O(1/{\rm dist}(\mu,{\rm Spec}(\Delta)))$ and the smoothness increasing theorems for solutions to elliptic equations (including the weighted estimates near vertices, see Chapters 3,4, \cite{NP}), formula (\ref{parametrix rhs estimate}) implies that
\begin{equation}
\label{parametrix estimate}
|\tilde{R}_\mu(x,x')|+|\partial_z\tilde{R}_\mu(x,x')|+|\partial_{z'}\tilde{R}_\mu(x,x')|=O(|\mu|^{-\infty}), \qquad \Re\mu\to-\infty
\end{equation}
uniformly in $x\in X$ and $x'\in\mathbb{K}_j(\epsilon)$; here $z,z'$ are arbitrary holomorphic coordinates of $x,x'$ on $X$, respectively. 

Substituting (\ref{parametrix resolvent}) into the first formula in (\ref{psi mu}) and taking into account (\ref{parametrix estimate}) and (\ref{HK Cone}), one obtains the (admitting differentiation and uniform in $x\in\mathbb{K}_j(\epsilon)$) asymptotics
\begin{align}
\nonumber
\psi_\mu(x)=\frac{\mu}{2\pi}\Big[K_0(d\sqrt{-\mu})-\int\limits_{0}^{+\infty}\int\limits_{\mathscr{C}}{\rm exp}\Big(\frac{\mathfrak{-r}^2}{4t}\Big)\Xi_{\beta_j}\frac{e^{\mu t}dtd\vartheta}{4i\beta_j t}\Big]_{x'=x}&+O(|\mu|^{-\infty})=\\ \nonumber
=\frac{\mu}{4\pi}\int\limits_{0}^{+\infty}\Big[{\rm exp}\Big(\frac{-d^2}{4t}\Big)-\int\limits_{\mathscr{C}}{\rm exp}\Big(\frac{\mathfrak{-r}^2}{4t}\Big)\Xi_{\beta_j}\frac{e^{\mu t}d\vartheta}{2i\beta_j}\Big]_{x'=x}\frac{dt}{t}&+O(|\mu|^{-\infty})=\\ 
\label{psi asym}
=a_\mu(|\xi_j(x)|,\beta_j)+O(|\mu|^{-\infty}), \qquad \Re&\mu\to-\infty,
\end{align}
where 
\begin{align*}
a_\mu(r,\beta):=\frac{-\mu}{8\pi i\beta}\int\limits_{\tilde{\mathscr{C}}}d\vartheta\,{\rm cot}\Big(\frac{\pi\vartheta}{\beta}\Big) \int\limits_{0}^{+\infty}{\rm exp}\Big(\mu t-\frac{r^2{\rm sin}^2(\vartheta/2)}{t}\Big)\frac{dt}{t}
\end{align*}
and the contour $\tilde{\mathscr{C}}$ is obtained from $\mathscr{C}$ by removing the small circle centered at $\vartheta=0$. Formulas (\ref{psi asym}) and (\ref{Fourier}) imply
\begin{align*}
{\rm p.v.}\int\limits_{\mathbb{K}_j(\epsilon)}&\psi_\mu\dot{\phi}dS={\rm p.v.}\int\limits_{\mathbb{K}_j(\epsilon)}a_\mu(|\xi_j(z)|,\beta_j)\dot{\phi}(x)dS(x)+O(|\mu|^{-\infty})=\\
=&\tilde{A}_j\Bigg[\beta_j\int\limits_{0}^{\epsilon}a_\mu(r,\beta_j){\rm log}r\,rdr\Bigg]+A_j\Bigg[\beta_j\int\limits_{0}^{\epsilon}a_\mu(r,\beta_j)rdr\Bigg]+O(|\mu|^{-\infty})+\\
&+{\rm p.v.}\int\limits_{0}^{\epsilon}dr\,a_\mu(r,\beta_j)r^{\frac{2\pi k}{\beta_j}+1}\int\limits_{0}^{\beta_j}\sum_{k=-1,1,2,\dots}\big[c_{k,j}e^{\frac{2\pi k i\varphi}{\beta_j}}+c^\dag_{k,j}e^{\frac{-2\pi k i\varphi}{\beta_j}}\big]d\varphi
\end{align*}
Here the last term equals zero since the function $a_\mu(|\xi_j(x)|,\beta_j)=a_\mu(r,\beta_j)$ is rotationally symmetric. The substitution of the last formula into (\ref{zeta calc 2}) yields
\begin{equation}
\label{zeta calc}
\partial_\tau[{\rm log}({\rm det}\Delta/A)]=\sum_{j=1}^M[\tilde{A}_j\beta_j\tilde{\mathfrak{I}}'(\beta_j)+A_j\mathfrak{I}(\beta_j)],
\end{equation}
where $\tilde{A}_j$ and $A_j$ are the coefficients in expansion (\ref{Fourier}) for $\dot{\phi}$ and
\begin{equation}
\label{integrals}
\tilde{\mathfrak{I}}'(\beta):=\Bigg[\int\limits_{0}^{\epsilon}a_\mu(r,\beta){\rm log}r\,rdr\Bigg]^{\infty}, \quad \mathfrak{I}(\beta):=\Bigg[\beta\int\limits_{0}^{\epsilon}a_\mu(r,\beta)rdr\Bigg]^{\infty}.
\end{equation}

\section{Expressions for $\tilde{\mathfrak{I}}$ and $\mathfrak{I}$}
\label{heat kernel sec}
The asymptotics
\begin{align}
\label{integrated asymptotic const term}
\beta\int\limits_{0}^{\epsilon}a_\mu(r,\beta)rdr=\mathfrak{I}(\beta)+O(|\mu|^{-\infty}),  \qquad \mathfrak{I}(\beta)=-\frac{1}{12}\Big(\frac{\beta}{2\pi}-\frac{2\pi}{\beta}\Big)
\end{align}
of the second integral in (\ref{integrals}) is derived in Proposition 1, \cite{KoKorotTau} from the explicit expressions for the heat kernel obtained in \cite{Dowker1}. The asymptotics 
\begin{equation}
\label{log int asym}
\int\limits_{0}^{\epsilon}a_\mu(r,\beta){\rm log}r\,rdr=-\frac{\mathfrak{I}(\beta)}{2\beta}{\rm log}(-\mu)+\tilde{\mathfrak{I}}'(\beta)+O(|\mu|^{-\infty})
\end{equation}
of the first integral in (\ref{integrals}) is derived in Subsection 4.2., \cite{KoKoPolyPoly}, where
\begin{equation}
\label{log integral at infty}
\tilde{\mathfrak{I}}(\beta)=\frac{\gamma-1-{\rm log}(\beta/2)}{12}\Big(\frac{\beta}{2\pi}+\frac{2\pi}{\beta}\Big)+\frac{5\beta}{48\pi}-\mathcal{H}\int\limits_{0}^{+\infty}{\rm coth}(\pi l)\,{\rm coth}\Big(\frac{\beta l}{2}\Big)\frac{dl}{8l}
\end{equation}
and $\mathcal{H}$ denotes the Hadamard regularization of the diverging integral (obtained by replacing the lower limit with $\varepsilon>0$, eliminating the term proportional $\varepsilon^{-1}$ and then passing to the limit as $\varepsilon\to+0$). 

Formulas (\ref{integrated asymptotic const term}), (\ref{log int asym}) and (\ref{log integral at infty}) show that zeta function (\ref{zeta of 2}) does indeed have the asymptotics of form (\ref{F asymp}); they also provide the expressions for coefficients (\ref{integrals}) in (\ref{zeta calc}).

\section{Formulas for coefficients $A_j$ and $\tilde{A}_j$}
\label{sec coeff}
\subsection{The case $\tau=z'_i$}
Let $\tau=z'_i$. Then (\ref{Fourier coef}), (\ref{u j explict}), and the equality $-2\partial_{z}\Im f(z)=\partial_{z}(if(z))$ for holomorphic $f$ yield $\tilde{A}_1=\dots=\tilde{A}_M=0$ and
\begin{align}
\label{A j ne i z i}
\begin{split}
A_j=[\partial_{z'_i}\Re u_j](z'_j)=-b_i\partial_{z'_i}&\Big({\rm log}\Big[\frac{E(z'_j,z'_i)}{E(p_1,z'_i)}\Big]+Q(z'_i,z'_j)\Big),
\end{split}
\end{align}
for $j\ne i$ where $Q$ is given by (\ref{Q form}). Similarly, (\ref{Fourier coef}), (\ref{u j explict}), and the asymptotics $E(z,z')=z'-z+O((z'-z)^3)$ of the prime form yield
\begin{align}
\label{A i z i}
\begin{split}
A_i=\partial_{z'_i}\Bigg[(b_i+2){\rm log}E(p_1,z_i)&+2{\rm log}\sigma(z'_i,p_0)-Q(z'_i,z'_i)+\\
+\sum_{j\ne i}b_j&\Big({\rm log}E(z'_i,z'_j)+Q(z'_i,z'_j)\Big)\Bigg]\frac{b_i}{b_i+2}.
\end{split}
\end{align}
Substituting (\ref{A j ne i z i}), (\ref{A i z i}) and (\ref{integrated asymptotic const term}) into (\ref{zeta calc}) and taking into account that
$$\frac{\beta_j}{2\pi}-\frac{2\pi}{\beta_j}=b_j\frac{b_j+2}{b_j+1}, \qquad \frac{b_j+2}{b_j+1}-\frac{b_i}{b_i+1}=2\pi\Big(\frac{1}{\beta_j}+\frac{1}{\beta_i}\Big)$$
one arrives at
\begin{align*}
\partial_{z'_i}[{\rm log}&({\rm det}\Delta/A)]=\sum_{j=1}^MA_j\mathfrak{I}(\beta_j)=\\
=\partial_{z'_i}&\Bigg[\frac{\pi}{6}\sum_{j\ne i} b_ib_j\Big(\frac{1}{\beta_i}+\frac{1}{\beta_j}\Big)\big({\rm log}E(z'_i,z'_j)+Q(z'_i,z'_j)\big)+\frac{\pi b^2_i}{6\beta_i}Q(z'_i,z'_i)+\\
&+\sum_{j=1}^M \mathfrak{J}(\beta_j)b_i{\rm log}E(p_1,z'_i)-\frac{2\pi b_i^2}{3\beta_i}{\rm log}|\sigma(z'_i,p_0)|\Bigg].
\end{align*}
Hence we have
\begin{align}
\label{Main formula 1}
\begin{split}
{\bm D}:=&\frac{\pi}{6}\Big[\sum_{i<j} b_ib_j\Big(\frac{1}{\beta_i}+\frac{1}{\beta_j}\Big){\rm log}|E(z'_i,z'_j)|^2+\frac{1}{2}\sum_{i,j} b_ib_j\Big(\frac{1}{\beta_i}+\frac{1}{\beta_j}\Big)Q(z'_i,z'_j)\Big]+\\
+&\sum_{i=1}^M\Big[\sum_{j=1}^M\mathfrak{J}(\beta_j)b_i{\rm log}|E(p_1,z'_i)|^2-\frac{2\pi b_i^2}{3\beta_i}{\rm log}|\sigma(z'_i,p_0)|\Big]={\rm log}\frac{{\rm det}\Delta}{A\mathbf{C}(\vec{\beta},C)}
\end{split}
\end{align}
where $Q$ is given by (\ref{Q form}) and $\mathbf{C}(\vec{\beta},C)$ is independent of $(z'_1,\dots,z'_M)^T=\vec{z}'$. Cumbersome but straightforward calculations based on formula (\ref{Troy 2}) show that ${\bm D}$ is independent of the choice of the local coordinates $z'_k$ representing $P_k$; thus, so is $\mathbf{C}(\vec{\beta},C)$.

\subsection{The case $\tau=C$}
Let $\tau=C$. Then (\ref{Fourier coef}) and (\ref{u j explict}) read $\tilde{A}_1=\dots=\tilde{A}_M=0$ and $A_1=\dots=A_M=-1/C$. Comparing (\ref{zeta calc}) with (\ref{Main formula 1}) and (\ref{zeta calc}) yields
\begin{equation}
\label{scaling formula}
\mathbf{C}(\vec{\beta},C)=\mathbf{C}(\vec{\beta})C^{-\sum_{j=1}^M \mathfrak{I}(\beta_j)},
\end{equation}
where $\mathbf{C}(\vec{\beta})$ is independent of $\vec{z}'$ and $C$ (provided that $z_0$ and the choice of the local coordinates $z'_k$ of vertices are fixed).

\subsection{The case $\tau=\beta_i$, $\,\,\beta_i+\beta_1={\rm const}$}
Let $\tau=\beta_i$ and $\beta_i+\beta_1={\rm const}$; the last constraint is needed to obey the Gauss-Bonnet formula $\sum_{j=1}^M b_j=2g-2$. Calculating $\partial_\tau[{\rm log}({\rm det}\Delta/A)]$ with help of (\ref{zeta calc 1}), (\ref{Fourier coef}) and (\ref{u j explict}) and comparing the result with equality (\ref{Main formula 1}) differentiated in $\tau$, one arrives at
\begin{equation}
\label{beta part DE}
\partial_{\beta_i}\Big({\rm log}\mathbf{C}(\vec{\beta})-2\tilde{\mathfrak{J}}(\beta_i)\Big)+\frac{2\mathfrak{J}(\beta_i)}{\beta_i}{\rm log}\Big(\frac{\beta_i}{2\pi}\Big)=0
\end{equation}
for $i=1,\dots,M$. Since all the terms depending on the positions $z'_k$ of vertices and the scaling factor $C$ should be canceled by each other during the above comparison, it can be considered as a (cumbersome but straightforward) cross-check of (\ref{Main formula 1}). The integration of (\ref{beta part DE}) yields
\begin{align}
\label{beta part}
\begin{split}
\mathbf{C}(\vec{\beta})&={\bf c}_0\prod_{j=1}^M{\bf d}^2(\beta_j),\\
{\bf d}(\beta)&={\rm exp}\Big(\tilde{\mathfrak{J}}(\beta)-\tilde{\mathfrak{J}}(2\pi)+\mathfrak{J}(\beta)+\frac{1}{12}\Big(\frac{\beta}{2\pi}+\frac{2\pi}{\beta}\Big){\rm log}\Big(\frac{\beta}{2\pi}\Big)\Big),
\end{split}
\end{align}
where ${\bf c}_0$ is independent of $\vec{z}'$, $C$, and $(\beta_1,\dots,\beta_M)^T=\vec{\beta}$ (provided that $z_0$ and the choice of the local coordinates $z'_k$ of vertices are fixed).

Combining formulas (\ref{Main formula 1}), (\ref{scaling formula}), and (\ref{beta part}), one arrives at
\begin{align}
\label{Main fomula 2}
\begin{split}
{\rm log}&\frac{{\rm det}\Delta}{A}={\bm D}+{\rm log}{\bf c}_0+\sum_{j=1}^M\Big[2{\rm log}{\bf d}(\beta_j)-\mathfrak{I}(\beta_j){\rm log}C\Big].
\end{split}
\end{align}

\subsection{The cases $\tau=p_0$ and $\tau=p_1$}
Differentiating (\ref{Main fomula 2}) with respect to $\tau=p_1$ and taking into account (\ref{zeta calc}), (\ref{Fourier coef}), (\ref{u j explict}), one arrives at 
$$\partial_{p_1}{\rm log}{\bf c}_0=0.$$ 
Similarly, let us differentiate (\ref{Main fomula 2}) with respect to $\tau=p_0$. Due to the multiplication law
$$\sigma(z,p_0)\sigma(p_0,p_1)=\sigma(z,p_1)=1/\sigma(p_1,z),$$
(see (1.13), \cite{Fay}), $\partial_{p_{0}}{\rm log}|\sigma(z,p_0)|$ is independent on $z$. Then,
\begin{align*}
\partial_{p_{0}}{\rm log}|\sigma(\cdot,p_0)|\,\Big[\sum_{i=1}^M\Big(\frac{\beta_i}{2\pi}-\frac{2\pi}{\beta_i}\Big)\Big]\frac{1}{3}=\partial_{p_0}[{\rm log}({\rm det}\Delta/A)]=\\
=\partial_{p_0}[{\bm D}+{\rm log}{\bf c}_0]=\partial_{p_0}\Big[\partial_{p_{0}}{\rm log}|\sigma(\cdot,p_0)|\,\sum_{i=1}^M\frac{b_i^2}{3(b_i+1)}+{\rm log}{\bf c}_0\Big].
\end{align*}
Thus, the ${\bf c}_0$ depends on the choice of $p_0$ and obeys the transformation law
$${\bf c}_0(p_0)=|\sigma(p_0,p'_0)|^{\frac{4(1-g)}{3}}{\bf c}_0(p'_0).$$

\section{Genus 1 case}
\label{sec tori}
Any polyhedral metric $m$ on the complex torus $X$ can be written in form (\ref{metric}), where
\begin{align}
\label{polyhe tori}
\phi(z)=-\sum_{k=1}^M 2b_k\Bigg[\frac{\pi\Im P_k \Im z}{\Im{\mathbb B}}+ {\rm log}\Big|\theta\left[^{1/2}_{1/2}\right](z-P_k)\Big|\Bigg]-{\rm log}C.
\end{align}
Here $\mathcal{D}=\sum_{k=1}^M b_k P_k$ is the generalized divisor of the metric, $\mathbb{B}$ is the $b$-period of $X$, $\theta[\chi]$ is the theta-function with characteristics $\chi$, $C$ is the scaling factor. Since the reasoning of Sections \ref{sec eig},\ref{sec det} is valid for any genus, the determinant of the Friedrichs Laplacian on $(X,m)$ obeys (\ref{zeta calc}). 

Let $\tau=P_i$ ($i=1,\dots,M$), then
\begin{align*}
A_j&=\partial_{P_i}\Big[\frac{\pi\Im P_i \Im P_j}{\Im{\mathbb B}}+{\rm log}\Big|\theta\left[^{1/2}_{1/2}\right](P_j-P_i)\Big|\Big](-2b_i) \quad \qquad (j\ne i),\\
A_i&=\partial_{P_i}\Bigg[\sum_{j\ne i}^M b_j\Big(\frac{\pi\Im P_j \Im P_i}{\Im{\mathbb B}}+{\rm log}\Big|\theta\left[^{1/2}_{1/2}\right](P_i-P_j)\Big|\Big)-\frac{\pi\Im P_i \Im P_i}{\Im{\mathbb B}}\Bigg]\frac{2b_i}{b_i+2}
\end{align*}
(here we used formulas (\ref{polyhe tori}) and (\ref{Fourier coef}), (\ref{u j explict}) and the fact that the theta functions with odd hall-integer characteristics are odd). The substitution of the above coeeficients into (\ref{zeta calc}) yields $[{\rm log}({\rm det}\Delta/A)-\mathbf{D}=\mathbf{C}(\vec{\beta},C)$, where 
\begin{align}
\label{torus D}
\begin{split}
\mathbf{D}=\frac{\pi}{3}\Bigg[\sum_{i<j} b_ib_j\Big(\frac{1}{\beta_i}+\frac{1}{\beta_j}\Big)&{\rm log}\Big|\theta\left[^{1/2}_{1/2}\right](P_i-P_j)\Big|+\\
+&\frac{1}{2}\sum_{i,j} b_ib_j\Big(\frac{1}{\beta_i}+\frac{1}{\beta_j}\Big)\frac{\pi\Im P_j \Im P_i}{\Im{\mathbb B}}\Bigg].
\end{split}
\end{align}
The same formulas with $\tau=C$ lead to (\ref{scaling formula}). Finally, let $\tau=\beta_i$ and $\beta_i+\beta_1={\rm const}$. Then formulas (\ref{zeta calc}), (\ref{polyhe tori}) and (\ref{Fourier coef}), (\ref{u j explict}) imply
$$\partial_{\beta_i}\Big({\rm log}\mathbf{C}(\vec{\beta})-2\tilde{\mathfrak{J}}(\beta_i)\Big)+\frac{2\mathfrak{J}(\beta_i)}{\beta_i}{\rm log}\Big(\frac{\beta_i}{2\pi}\Big|\theta'\left[^{1/2}_{1/2}\right](0)\Big|\Big)=0,$$
whence $\mathbf{C}(\vec{\beta})={\bf c}_0(\mathbb{B})\prod_{j=1}^M{\bf d}^2(\beta_j)$, where ${\bf c}_0(\mathbb{B})$ depends only on the moduli and
\begin{align}
\label{beta part torus}
\begin{split}
{\bf d}(\beta)={\rm exp}\Bigg(\tilde{\mathfrak{J}}(\beta)-\tilde{\mathfrak{J}}(2\pi)+\mathfrak{J}(\beta)+\frac{1}{12}\Big(\frac{\beta}{2\pi}+\frac{2\pi}{\beta}\Big){\rm log}\Big(\frac{\beta}{2\pi}\Big)+\\
+\frac{1}{12}\Big(\frac{\beta}{2\pi}+\frac{2\pi}{\beta}-2\Big){\rm log}\Big(\Big|\theta'\left[^{1/2}_{1/2}\right](0)\Big|\Big)\Bigg).
\end{split}
\end{align}
Here the additional (independent of $\beta$) terms $\tilde{\mathfrak{J}}(2\pi)$ and $-2$ are added manually to obey ${\bf d}(2\pi)=0$. Thus, we have proved formula (\ref{Main fomula 2}) with $\mathbf{D}$ and ${\bf d}$ given by (\ref{torus D}) and (\ref{beta part torus}), respectively. Finally, the substitution $\beta_1=\dots=\beta_M$, $C=1$ into (\ref{polyhe tori}) and (\ref{Main fomula 2}) leads to the equality 
$${\bf c}_0(\mathbb{B})={\rm det}\Delta^{(|dz|^2)}/A^{(|dz|^2)},$$
where the Laplacian $\Delta^{(|dz|^2)}$ and the area $A^{(|dz|^2)}$ correspond to the metric $m=|dz|^2$ on $X$. In view of the Ray-Singer formula \cite{RS}
\begin{equation*}
\frac{{\rm det}\Delta^{(|dz|^2)}}{A^{(|dz|^2)}}=\Im\mathbb{B}\, |\eta(\mathbb{B})|^4
\end{equation*}
(where $\eta$ is the Dedekind eta-function), one gets
\begin{equation}
\label{torical const}
{\bf c}_0(\mathbb{B})=\Im\mathbb{B}\, |\eta(\mathbb{B})|^4.
\end{equation}


\begin{thebibliography}{99}

\bibitem{BFK} 
	\newblock{Burghelea, D.; Friedlander, L.; Kappeler, T.}
	\newblock{\em Mayer-Vietoris type formula for determinants of elliptic differential operators.}
	\newblock{\em J. Funct. Anal.} 107 (1992), no. 1, 34–65.
	
	
	\bibitem{Carslaw}
	\newblock{ H.S. Carslaw.}
	\newblock {\em The Green's Function for a Wedge of any Angle, and Other Problems in the Conduction of Heat}.
	\newblock {\em Proc. Lond. Math. Soc.} 2(8) (1910): 365--374. DOI: \url{https://doi.org/10.1112/plms/s2-8.1.365}
	
	\bibitem{Dowker1}
	\newblock{J.S. Dowker.}
	\newblock {\em Quantum field theory on a cone.}
	\newblock {\em J. Phys. A: Math.} 10(1) (1977): 115--124. DOI: \url{https://doi.org/10.1088/0305-4470/10/1/023}
	
	\bibitem{Dowker}
	\newblock{J. S. Dowker.}
	\newblock{\em Effective action in spherical domains.} 
	\newblock{\em 
Commun. Math. Ph.} 162 (1994) 633-647
	
	
	\bibitem{Fay}
	\newblock{John Fay.}
	\newblock {\em Kernel functions, analytic torsion, and moduli spaces.}
	\newblock {\em Memoirs of the AMS} 464, Providence, Rhode Island (1992), 123 p. ISBN: 082182550X.
	
	
	\bibitem{KokLN}
	Alexey Kokotov.
	\newblock {On the Spectral Theory of the Laplacian on Compact Polyhedral Surfaces of Arbitrary Genus. In: Computational Approach to Riemann Surfaces.}
	\newblock {\em Lecture Notes in Mathematics} (2013): 227-253. 		\url{https://doi.org/10.1007/978-3-642-17413-1_8}
	
	\bibitem{KokPAMS}
	\newblock{A. Kokotov.}
	\newblock{\em Polyhedral surfaces and determinant of Laplacian.}
	\newblock{Proceedings of AMS.} Vol. 141, n. 2, 2013, p. 725-735
	
	
	\bibitem{KoKo-poly}
	\newblock{A. Kokotov, D. Korikov.} 
	\newblock{\em Regularized $\zeta_{\Delta}(1)$ for polyhedra.}
	\newblock{arXiv:2502.03351}	
	
	\bibitem{KokKorWZ}
	\newblock{A. Kokotov, D. Korikov.} 
	\newblock{\em On a polygon version of Wiegmann-Zabrodin formula.}
	\newblock{arXiv:2503.13718}
	
	\bibitem{KoKoPolyPoly}
	\newblock{A. Kokotov, D. Korikov.} 
	\newblock{\em On an infinitesimal Polyakov formula for genus zero polyhedra.}
	\newblock{arXiv:2504.17652}
	   
	
	\bibitem{KoKorotTau}
	\newblock{A. Kokotov, D. Korotkin.}
	\newblock {\em Tau-functions on spaces of Abelian differentials and higher genus generalizations of Ray-Singer formula}.
	\newblock {\em J. Differential Geom.} 82 (2004), 35--100.
	
	\bibitem{MNP}
	\newblock{V. Mazy'a, S. Nazarov, B. Plamenevskii.}
	\newblock{\em Asymptotic Theory of Elliptic Boundary Value Problems in Singularly Perturbed Domains.} Springer, 2000
	
	
	\bibitem{NP}
	\newblock{S. Nazarov, B. Plamenevskii.}
	\newblock{\em lliptic Problems in Domains with Piecewise Smooth Boundaries.} Berlin, New York: De Gruyter, 1994. 
	
	\bibitem{RS}
	\newblock{D. B. Ray, I. M. Singer.}
	\newblock {Analytic torsion for complex manifolds}.
	\newblock {\em Ann. of Math.} 98(1) (1973): 154--177. 
	
	\bibitem{Spreafico}
	\newblock{M. Spreafico.}
	\newblock{\em Zeta function and regularized determinant on a disk and on a cone.}
	\newblock{Journal of Geometry and Physics.}  54 (2005) 355-371
	
	\bibitem{Troy} 
	\newblock{Marc Troyanov.}
	\newblock {Les surfaces euclidiennes \`a singularit\'es coniques}.
	\newblock {\em L'Enseignement Math\'ematique} 32(2) (1986): 79--94. DOI: \url{https://doi.org/10.5169/seals-55079}

\end{thebibliography}
\end{document}